\documentclass[11pt, reqno]{amsart}
\usepackage{amssymb, amsthm, amsmath, amsfonts}
\usepackage[alphabetic]{amsrefs}
\usepackage{verbatim}

\usepackage{graphics}

\newcommand{\demo}{\begin{proof}}
\newcommand{\edemo}{\end{proof}}
\newcommand{\demoname}[1]{\begin{proof}[#1]}
\newcommand{\edemoname}{\end{proof}}

\newcommand{\stepname}{Step}

\newcommand{\setstepn}[1]{\setcounter{stepn}{#1}} 

\theoremstyle{plain}
\newtheorem{theorem}{Theorem}[section]

\newtheorem{corollary}[theorem]{Corollary}
\newtheorem{lemma}[theorem]{Lemma}

\theoremstyle{definition}
\newtheorem{example}[theorem]{Example}
\newtheorem*{definition}{Definition}
\newtheorem{step}{Step}
\newtheorem{casee}{Case}
\newtheorem{stepn}{\stepname}
\newtheorem*{stepnn}{\stepname}

\newcommand{\thm}{\begin{theorem}}
\newcommand{\ethm}{\end{theorem}}

\newcommand{\expl}{\begin{example}}
\newcommand{\eexpl}{\qex\end{example}}
\newcommand{\defn}{\begin{definition}}
\newcommand{\edefn}{\qef\end{definition}}
\newcommand{\stp}{\begin{step}}
\newcommand{\estp}{\end{step}}
\newcommand{\cse}{\begin{casee}}
\newcommand{\ecse}{\end{casee}}
\newcommand{\stpn}[1]{\renewcommand{\stepname}{#1}\begin{stepn}}
\newcommand{\estpn}{\end{stepn}}
\newcommand{\stpnn}[1]{\renewcommand{\stepname}{#1}\begin{stepnn}}
\newcommand{\estpnn}{\end{stepnn}}

\newcommand{\coro}{\begin{corollary}}
\newcommand{\ecoro}{\end{corollary}}
\newcommand{\lem}{\begin{lemma}}
\newcommand{\elem}{\end{lemma}}

\providecommand{\qexsymbol}{$\lozenge$}%
\newcommand{\mathqex}{\quad\hbox{\qexsymbol}}
\DeclareRobustCommand{\qex}{%
  \ifmmode \mathqex
  \else
    \leavevmode\unskip\penalty9999 \hbox{}\nobreak\hfill
    \quad\hbox{\qexsymbol}%
  \fi
}
\providecommand{\qefsymbol}{$\triangle$}%
\newcommand{\mathqef}{\quad\hbox{\qefsymbol}}
\DeclareRobustCommand{\qef}{%
  \ifmmode \mathqef
  \else
    \leavevmode\unskip\penalty9999 \hbox{}\nobreak\hfill
    \quad\hbox{\qefsymbol}%
  \fi
}

\newcommand{\enum}{\begin{enumerate}}
\newcommand{\eenum}{\end{enumerate}}

\newcommand{\nn}{\mathbb{N}}
\newcommand{\zz}{\mathbb{Z}}
\newcommand{\rrr}[1]{{\mathbb{R}}^{#1}}
\newcommand{\rr}{\mathbb{R}}

\newcommand{\nd}[1]{$#1$\nobreakdash-\hspace{0pt}}

\newcommand{\kd}[1]{d({#1})}
\newcommand{\kdk}[2]{\kd {{#1}, {#2}}}

\newcommand{\kb}[1]{S({#1})}
\newcommand{\kbk}[2]{\kb {{#1}, {#2}}}

\newcommand{\stcurv}{standard curvature}

\newcommand{\sdcurv}{stratified curvature}

\newcommand{\modsdcurv}{modified stratified curvature}

\newcommand{\sddef}{generalized angle defect}

\newcommand{\seccc}{stratified Euler characteristic}

\newcommand{\aang}{\alpha}
\newcommand{\angg}[2]{\aang({#1}, {#2})}
\newcommand{\angge}[2]{\aang^*({#1}, {#2})}

\newcommand{\simpn}[1]{\nd{{#1}}dimensional simplicial complex}

\DeclareMathOperator{\starpl}{star}
\DeclareMathOperator{\linkpl}{link}

\newcommand{\stark}[2]{\starpl\,({#1}, {#2})}
\newcommand{\linkk}[2]{\linkpl\,({#1}, {#2})}

\newcommand{\svf}{vertex-supported function}
\newcommand{\svfs}{vertex-supported functions}

\newcommand{\scvf}{simplicial-complex-supported function}
\newcommand{\scvfs}{simplicial-complex-supported functions}

\newcommand{\sst}{simplicial surface}
\newcommand{\ssts}{simplicial surfaces}

\newcommand{\gld}{geometrically locally determined}
\newcommand{\scl}{star-closed}

\newcommand{\simpt}{\nd{2}dimensional simplicial complex}
\newcommand{\simpts}{\nd{2}dimensional simplicial complexes}

\newcommand{\phiw}[2]{\phi({#1}, {#2})}
\newcommand{\psiw}[2]{\psi({#1}, {#2})}

\newcommand{\muw}[2]{\mu({#1}, {#2})}
\newcommand{\phia}[1]{\phi({#1})}

\newcommand{\Lambdaa}[1]{\Lambda({#1})}
\newcommand{\ip}[2]{\langle {#1}, {#2}\rangle}
\newcommand{\linkord}[2]{\mathcal{O}({#2}, {#1})}
\newcommand{\flaplinkord}[2]{\nflap {\linkord {#1}{#2}}}

\newcommand{\verts}[1]{{#1}^{(0)}}

\newcommand{\cad}{classical angle defect}

\newcommand{\dgbt}{Descartes-Gauss-Bonnet Theorem}

\newcommand{\tsimp}{\nd{2}simplex}
\newcommand{\tsimps}{\nd{2}simplices}

\newcommand{\nflap}[1]{\nd{#1}flap}
\newcommand{\nflaps}[1]{\nd{#1}flaps}
\newcommand{\evert}{end-vertex}
\newcommand{\everts}{end-vertices}
\newcommand{\ord}[2]{\mathsf{e}({#1}, {#2})}
\newcommand{\mcals}{\mathcal{S}}
\newcommand{\mcalt}{\mathcal{T}}

\begin{document}

\title{A Characterization of the Angle Defect and the Euler Characteristic in Dimension $2$\\
Preliminary Draft}
\author{Ethan D.\ Bloch}
\address{Bard College\\
Annandale-on-Hudson, NY 12504\\
U.S.A.}
\email{bloch@bard.edu}
\date{}
\subjclass[2000]{Primary 52B70}
\keywords{angle defect, Euler characteristic, simplicial complex, curvature, Descartes-Gauss-Bonnet Theorem}
\begin{abstract}
The angle defect, which is the standard way to measure curvature at the vertices of polyhedral surfaces, goes back at least as far as Descartes.  Although the angle defect has been widely studied, there does not appear to be in the literature an axiomatic characterization of the angle defect.  We give a characterization of the angle defect for \ssts, and we show that variants of the same characterization work for two known approaches to generalizing the angle defect to arbitrary \simpts.  Simultaneously, we give a characterization of the Euler characteristic on \simpts\ in terms of being \gld.
\end{abstract}
\maketitle

\markright{A CHARACTERIZATION OF THE ANGLE DEFECT IN DIMENSION $2$}

\section{Introduction}
\label{secINTR}

This paper is concerned with two related questions.  The first issue concerns the curvature at a vertex $v$ of a triangulated polyhedral surface $M$, which is given by $\kdk vM = 2\pi - \sum_{\alpha \ni v} \alpha$, where the $\alpha$ are the angles at $v$ of the triangles containing $v$.  This curvature function, which we refer to as the ``{\cad},''
goes back at least as far as Descartes (see \cite{FE}).  The \cad\ satisfies various properties one would expect a curvature function on polyhedra to satisfy.  For example, the angle defect is locally defined; it is invariant under simplicial local isometries (that is, functions that preserve the lengths of edges); it is zero at a vertex that has a flat star; it is invariant under subdivision; and it satisfies the polyhedral \dgbt, which says $\sum_v \kdk vM = 2\pi \chi(M)$, where the summation is over all the vertices of $M$, and where $\chi(M)$ is the Euler characteristic of $M$.  (We refer to this theorem as the ``\dgbt'' rather than just ``Descartes' Theorem'' because Descartes' version was for convex polytopes only, and did not explicitly involve the Euler characteristic for \nd{2}dimensional polyhedra; there appears to be some dispute in the literature as to whether or not Descartes was implicitly aware of the Euler characteristic.)

The angle defect (also known as the angle deficiency), and related constructs involving sums of angles in polyhedra, have been widely studied, both in the classical situation, as well as in higher dimensions.  It has been studied in the case of convex polytopes from a combinatorial approach, for example in \cite{SH2} and \cite{GR2}; more generally, for the wider study of angle sums in convex polytopes, see for example \cite{GR1}*{Chapter~14}, \cite{SH}, \cite{P-S} and \cite{MC}.  In \cite{G-S2} a \dgbt\ is proved for the angle defect in polytopes with underlying spaces that are manifolds.  This approach has been generalized to arbitrary simplicial complexes in \cite{BL4}, and further studied in \cite{BL10} and \cite{BL11}.  A different approach to generalizing the \cad\ has been studied extensively from a differential geometric point of view; see, among others, \cite{BA1}, \cite{WINT}, \cite{BUDA}, \cite{CHEE}, \cite{C-M-S} and \cite{ZA}.  One treatment of curvature of polyhedra that has some of the advantages of all the approaches cited above is in \cite{BA3}, which uses curvatures functions based on critical points (similarly to \cite{BA1}), but this time using projection maps $\rrr{n} \to \rrr{m}$, which leads to curvature functions related to the Grassman angles of \cite{GR2}, and which are located at all simplices, and which directly generalizes \stcurv; moreover, an angle defect type formula for curvature is obtained using projection maps $\rrr{n} \to \rrr{n-1}$.  (The angle defect and its generalizations treated in the above references, and which we discuss in this note, are geometric in nature, depending upon the measurement of angles; we will not be discussing the combinatorial approach to curvature of simplicial complexes found in \cite{FORM6}.) 

Although the angle defect has been widely studied, there does not appear to be in the literature an axiomatic characterization of the angle defect.  Such a characterization would be useful not only for gaining insight into the angle defect, but also to help distinguish between different generalizations of the \cad\ to arbitrary polyhedra.  In the present paper we give a characterization of the \cad\ for \ssts, and we show that variants of the same characterization work for the two approaches to generalizing the \cad\ to arbitrary \simpts\ , as found in \cite{BL4} on the one hand, and in \cite{BA1} et al.\ on the other.  

The second issue we discuss concerns the Euler characteristic, which is intimately connected to the \cad\ by the \dgbt.  In \cite{FORM5}, which follows \cite{LEVI}, it is shown that the Euler characteristic is the unique locally determined numerical invariant  of finite simplicial complexes that assigns the same number to every cone, where in this context simplicial complexes are considered to be the same if they are combinatorially equivalent.  A real-valued function $\rho$ on the set of all finite simplicial complexes is ``locally determined'' if there is another real-valued function $h$ on the set of all finite simplicial complexes such that for each simplicial complex $K$, we have $\rho (K) = \sum_v h(\linkk vK)$, where $\linkk vK$ denotes the link of $v$ in $K$, and where the summation is over all the vertices of $K$.  This last condition certainly has some resemblance to the \dgbt, but there is a substantial difference between the results of \cite{LEVI} and \cite{FORM5} on the one hand, and the \dgbt\ on the other: the nature of being locally determined in \cite{LEVI} and \cite{FORM5} is purely combinatorial, whereas in the \dgbt\ the Euler characteristic is locally determined by a geometric quantity, namely the \cad, which depends not only on the combinatorial nature of the link of each vertex, but on the geometry of the embedding of the star of each vertex.  One might therefore think of the results of \cite{LEVI} and \cite{FORM5} as characterizing the Euler characteristic among those functions that are ``combinatorially locally determined.''  In the present paper we will given an analogous characterization of the Euler characteristic in the \nd{2}dimensional case among those functions that are ``\gld,'' a term that will be defined below.  Not surprisingly, our characterizations of both the angle defect and the Euler characteristic in dimension $2$ are simply different ways of looking at the same result.

\section{Statement of Results}
\label{secSTAT}

We start with some assumptions, notation and definitions.  Throughout this paper we will restrict our attention to simplicial complexes, all of which are \nd{2}dimensional, finite, and are embedded in Euclidean space (which we will not name when it is not necessary).  Different embeddings of combinatorially equivalent simplicial complexes will be considered as different simplicial complexes (in contrast to \cite{LEVI} and \cite{FORM5}, whose approach is combinatorial rather than geometric).  We use the term simplicial surface to mean a \nd{2}dimensional simplicial complex the underlying space of which is a compact surface without boundary.

Let $K$ be a simplicial complex.  We let $|K|$ denote the underlying space of $K$.  We write $\verts K$ to denote the set of vertices of $K$, and $f_i(K)$ to denote the number of \nd{i}simplices of $K$ for each possible value of $i$.  If $v, w \in \verts K$, and if there is a \nd{1}simplex of $K$ with vertices $v$ and $w$, we let $\ip vw$ denote this \nd{1}simplex, and we let $\linkord wv$ denote the number of \nd{2}simplices of $K$ that contain $\ip vw$.  If $\sigma \in K$, we use $\linkk {\sigma}K$ and $\stark {\sigma}K$ to denote the link and star of $\sigma$ in $K$ respectively.  (For basic definitions in PL topology, see for example \cite{GL}*{vol.\ I} and \cite{HU}, although  the latter uses the notation $\overline{\textrm{star}}$ instead of $\textrm{star}$.)

For the sake of convenience, we adopt the convention that all angles in \tsimps\ are normalized so that the circumference of the unit circle is $1$.  For any \tsimp\ $\sigma^2$ in Euclidean space, and any vertex $v$ of $\sigma^2$, we let $\angg v{\sigma^2}$ denote the (interior) angle in $\sigma^2$ at $v$, where by normalization such an angle is always a number in $[0, \frac 12]$.  Hence the $2\pi$ will drop out of our statement of the \dgbt.  

The following definition gives the two general types of functions of which the Euler characteristic and the \cad\ are examples, respectively.

\defn 
Let $\mcalt$ be a set of \simpts.  A \textbf{\scvf} on $\mcalt$ is a function $\Lambda$ that assigns to every \simpt\ $K \in \mcalt$ a real number $\Lambdaa K$.  A \textbf{\svf} on $\mcalt$ is a function $\phi$ that assigns to every \simpt\ $K \in \mcalt$, and every vertex $v \in \verts K$, a real number $\phiw vK$.
\edefn

We now consider various properties of \svfs.  These properties, defined below, are all satisfied by the \cad, as can be verified easily.  The first property involves subdivision.  

\defn
Let $\mcalt$ be a set of \simpts, and let $\phi$ be a \svf\ on $\mcalt$.  We say that $\phi$ is \textbf{invariant under subdivision} if the following condition holds.  Let $K, J \in \mcalt$, and let $v \in \verts K$.  If $J$ is a subdivision of $K$, then $\phiw vK = \phiw vJ$.
\edefn

For our next property, which involves isometries, we need the following terminology.  Let $K$ and $L$ be \simpts.  We say that  $K$ and $L$ are simplicially isometric if there is a simplicial homeomorphism $|K| \to |L|$ that preserves the lengths of edges; such a map is called a simplicial isometry.  Moreover, let $v \in \verts K$ and let $w \in \verts L$.  We say that $\stark vK$ and $\stark wL$ are simplicially isometric if there is a simplcial isometry $|\stark vK| \to |\stark wL|$ that takes $v$ to $w$; any simplicial isometry between $\stark vK$ and $\stark wL$ will always be assumed to take $v$ to $w$.

\defn
Let $\mcalt$ be a set of \simpts, and let $\phi$ be a \svf\ on $\mcalt$.  We say that $\phi$ is \textbf{invariant under simplicial isometries of stars} if the following condition holds.  Let $K, L \in \mcalt$, let $v \in \verts K$ and let $w \in \verts L$.  If $\stark vK$ and $\stark wL$ are simplicially isometric, then $\phiw vK = \phiw wL$.
\edefn

Our next property involves continuity.  Suppose that $K$ and $\displaystyle \{K_n\}_{n=1}^\infty$ are combinatorially equivalent \simpts, and all are embedded in the same Euclidean space.  We can think of all these \simpts\ as embeddings of the same abstract simplicial complex.  We write $\lim_{n\to\infty} K_n = K$ to denote pointwise convergence of these embeddings; it suffices to verify convergence at the vertices of the abstract simplicial complex.

\defn
Let $\mcalt$ be a set of \simpts, and let $\phi$ be a \svf\ on $\mcalt$.  We say that $\phi$ is \textbf{continuous} if the following condition holds.  Let $\displaystyle \{K_n\}_{n=1}^\infty$ and $K$ be combinatorially equivalent \simpts\ in $\mcalt$, all embedded in the same Euclidean space.  Suppose $\lim_{n\to\infty} K_n = K$.  If $v \in \verts K$, and if we label the corresponding vertices of the $K_n$ as $v_n$, then $\lim_{n\to\infty}\phiw {v_n}{K_n} = \phiw vK$.
\edefn

Our final property is the analog of the \dgbt.

\defn
Let $\mcalt$ be a set of \simpts, let $\phi$ be a \svf\ on $\mcalt$, and let $\Lambda$ be a \scvf\ on $\mcalt$.  We say that $\phi$ \textbf{satisfies the \dgbt\ with respect to} $\Lambda$ if the following condition holds.  If $K \in \mcalt$, then $\sum_{v \in K} \phiw vK = \Lambdaa K$, where the sum is over all the vertices of $K$.  
\edefn

We can now state our results, which take place in two contexts, namely \ssts\ on the one hand, and the set of all finite \simpts\ on the other hand.  We first consider things from the point of view of the the \cad\ and its generalizations; shortly we will turn to the point of view of the Euler characteristic and its generalizations.  For \ssts, we have the following characterization of the \cad.

\thm\label{thmSURF} 
Let $\phi$ be a \svf\ on the set of all \ssts.  Then $\phi$ is the
\cad\ iff $\phi$ is invariant under simplicial isometries of stars and under subdivision, is continuous, and satisfies the \dgbt\ with respect to the Euler characteristic.
\ethm

It would be interesting to know whether the four hypotheses in Theorem~\ref{thmSURF} are all needed.  We do not have a complete answer to this question, though the following examples partially answer this question.

\expl\label{explH2}
\stpnn{(1)}\label{stpnnFRED}
It is clear that the \dgbt\ cannot be dropped in Theorem~\ref{thmSURF}, because the constantly zero \svf\ on the set of all \ssts\ satisfies the three other criteria of the theorem.
\estpnn

\stpnn{(2)}
If $K$ is a \sst, and if $v \in \verts K$, we let $\ord vK$ denote the number of edges of $K$ that contain $v$.  Define the \svf\ $\psi$ on the set of all \ssts\ by letting $\psiw vK = 1 - \frac 13\ord vK$ for any \sst\ $K$ and any vertex $v$ of $K$.  It is evident that $\psi$ is invariant under simplicial isometries of stars and is continuous, and it can also be verified that $\psi$ satisfies the \dgbt\ with respect to the Euler characteristic (using the fact that $\sum_{v \in K} \ord vK = 2f_1(K) = 3f_2(K)$, where the summation is over all vertices of $K$).  However, it is evident that $\psi$ is not invariant under subdivision, and therefore the invariance under subdivision criteria cannot be dropped from Theorem~\ref{thmSURF}.
\estpnn

\stpnn{(3)}
If $K$ is a \sst, we let $n(K)$ denote the number of vertices $v$ of $K$ for which the sum of the angles at $v$ is not equal to $1$; observe that $n(K)$ is never zero, because extreme vertices have positive \cad.  Define the \svf\ $\mu$ on the set of all \ssts\ by letting 
$$
\muw vK = 
\begin{cases}
\chi(K)/n(K),&\text{if the sum of the angles at $v$ is not $1$,}\\
0,&\text{otherwise,}
\end{cases}
$$
 for any \sst\ $K$ and any vertex $v$ of $K$.  It is clear that $\mu$ satisfies the \dgbt\ with respect to the Euler characteristic, and it can also be seen that $\mu$ is invariant under subdivision (because any new vertex in a subdivision of a \simpt\ $K$ is in the relative interior of a \nd{1}simplex or a \nd{2}simplex of $K$, and hence has angle sum equal to $1$).  We leave it to the reader to verify that $\mu$ is not invariant under simplicial isometries of stars and is not continuous.  Hence, these last two criteria cannot both be dropped from Theorem~\ref{thmSURF}.
\estpnn
\eexpl

Theorem~\ref{thmSURF} holds unchanged for the class of all \nd{2}dimensional simplicial pseudomanifolds (without boundary), and we omit the details.  The situation becomes more interesting when we go beyond pseudomanifolds, and look at the class of all \simpts, because for non-pseudomanifolds, there are (at least) two generalizations of the \cad, both of which are equal to the \cad\ when restricted to pseudomanifolds.

Both of these generalizations of the \cad\ work in all dimensions.  The first of these generalizations, which I refer to as ``\stcurv,'' is discussed, among others, in \cite{BA1}, \cite{WINT}, \cite{BUDA}, \cite{CHEE}, \cite{C-M-S} and \cite{ZA}.  It is very simple to define (though it does not directly resemble the \cad), and it satisfies our four properties.  In all dimensions, this type of curvature is concentrated at the vertices of simplicial complexes.  

The second generalization of the \cad, which I refer to as the ``\sddef,'' was defined in \cite{BL4}, and further studied in \cite{BL10} and \cite{BL11}.  This type of curvature, which more closely resembles the \cad\ than \stcurv, is a generalization of the higher dimensional angle defect for convex polytopes and manifolds studied, among others, in \cite{SH2}, \cite{GR2} and \cite{G-S2}.  In dimensions higher than $2$, the \sddef\ is not concentrated at the vertices of simplicial complexes, but rather is defined for each simplex of codimension at least $2$.  (A word on our terminology.  In order to compare our approach with \stcurv, we originally somewhat artificially concentrated our curvature at the vertices in Section~3 of \cite{BL4}, and called it ``\sdcurv.''  In Section~4 of \cite{BL4} we took the more natural approach that we are using at present, and referred to this approach by the unfortunate name ``\modsdcurv,'' which misses the point that in this approach we are really working with a pure angle defect.  Hence, in the present paper, we will use the better name ``\sddef,'' which we also use in \cite{BL10} and \cite{BL11}.  A detailed comparison of \stcurv\ with both \sdcurv\ and the \sddef\ may be found in \cite{BL4}*{Section~4}.)  

It turns out that in the present paper we will not ever need the actual definitions of \stcurv\ and the \sddef\ in the present paper---all we need is their properties.  Like \stcurv, the \sddef\ is invariant under simplicial isometries of stars and under subdivision, and it is continuous.  It also satisfies a \dgbt, though not with respect to the Euler characteristic, but with respect to a variant of the Euler characteristic, called the \seccc.  See \cite{BL4} for the definition of the \seccc, and a proof of the \dgbt\ for the \sddef.

The following theorem characterizes both these types of curvatures on the set of all \simpts.

\thm\label{thmSTCU} 
Let $\phi$ be a \svf\ on the set of all \simpts.  Then $\phi$ is \stcurv\ (respectively the \sddef) iff $\phi$ is invariant under simplicial isometries of stars and under subdivision, is continuous, and satisfies the \dgbt\ with respect to the Euler characteristic (respectively the \seccc).
\ethm

Theorem~\ref{thmSTCU} sheds light on how similar \stcurv\ and the \sddef\ are in dimension $2$.  These two types of curvature are less similar in higher dimensions, because \stcurv\ is concentrated at vertices and the \sddef\ is not.  A better understanding of the difference between these two types of curvature awaits characterization of \stcurv\ and the \sddef\ in higher dimensions.  Unfortunately, our proof of Theorem~\ref{thmSTCU} (which is really a corollary to Theorem~\ref{thmMAINA} stated below) does not generalize beyond the \nd{2}dimensional case.  It would be interesting to know whether the higher dimensional analogs of our results are nonetheless true, using a different method of proof.

Also, we note that Theorem~\ref{thmSTCU} does not imply Theorem~\ref{thmSURF}, because more is being assumed about $\phi$ in Theorem~\ref{thmSTCU} than in Theorem~\ref{thmSURF}.

We now turn our attention to the point of view of the Euler characteristic and its generalizations.  The following definition gives our geometric analog of the notion of being locally determined as discussed in \cite{LEVI} and \cite{FORM5}.

\defn
Let $\mcalt$ be a set of \simpts, and let $\Lambda$ be a \scvf\ on $\mcalt$.  We say that $\Lambda$ is \textbf{\gld} if there is a \svf\ $\phi$ on $\mcalt$ such that $\phi$ is invariant under simplicial isometries of stars and under subdivision, is continuous, and satisfies the \dgbt\ with respect to $\Lambda$.  If we need to specify $\phi$, we will say that $\Lambda$ is \gld\ by $\phi$.
\edefn

It is reasonable to expect that not every arbitrary \scvf\ will be \gld, and hence we should restrict our attention to those \scvfs\ that are well-behaved in some appropriate way.  In \cite{FORM5}, as seen in the title of that paper, the condition of being constant on the set of cones is used; this condition is weaker than the condition of being a homotopy invariant, which is used in \cite{LEVI}.  Because the \seccc\ of \cite{BL4}*{Section~2} is not a homotopy invariant (though it is a homeomorphism invariant), we adopt the approach of \cite{FORM5}, and will consider \scvfs\ that are constant on a number of different sets of \simpts, as we now discuss.

We will use the following standard terminology.  Let $D$ be a simplicial disk in $\rrr{2}$.  A pyramid on $D$ is the \sst\ obtained by coning on $D$ from a point in $\rrr{3}$ (called the apex of the pyramid) that is not in $\rrr{2}$, and then taking the boundary.  Let $R$ be a polygonal disk in $\rrr{2}$ (not necessarily subdivided into simplices).  A bipyramid on $R$ is the \sst\ obtained by suspending $D$ from two points in $\rrr{3}$ (called the apices of the bipyramid) that are not in $\rrr{2}$, and then taking the boundary.  We also need the following terms.

\defn
A \textbf{planar fan} is a simplicial disk in $\rrr{2}$ with no interior vertices that is triangulated as a cone from one of its boundary vertices.  Let $n \in \zz$ be such that $n \ge 0$.  An \textbf{\nflap n} with \textbf{\everts} $v$ and $w$ is a \simpt\ $K$ containing vertices $v$ and $w$ such that the following conditions hold: (1) $\ip vw$ is an edge of $K$; (2) $K$ has precisely $n$ \tsimps, each of which has $\ip vw$ as an edge; and (3) $K$ has no edges other than the \nd{1}faces of these $n$ \tsimps.
\edefn

Observe that there can be \nflaps n\ for any non-negative integer $n$.  A \nflap 0\ is a single edge, and a \nflap 1\ is a single \tsimp.  It is straightforward to verify that a \simpn 2\ $L$ is an \nflap n\ with \everts\ $v$ and $w$ iff $L = \stark vL = \stark wL$.  Also, if $K$ is a \simpt\ and $v \in \verts K$, and if $N = \stark vK$, then for any vertex $x$ of $\linkk vK$, it is seen that $\stark xN$ is a \flaplinkord xv.

As mentioned above, we will consider \scvfs\ that are constant on various sets of \simpts.  In particular, if $\Lambda$ is a \scvf, we will consider the case where $\Lambda$ is constant on the set of all planar fans, and we will write $\Lambdaa {\text{fan}}$; the case where $\Lambda$ is constant on the set of all \nflaps n, and we will write $\Lambdaa {\text{\nflap n}}$, for each $n \in \nn$; the case where $\Lambda$ is constant on the set of all pyramids and bipyramids, and we will write $\Lambdaa {\text{pyramid}}$; and the case where $\Lambda$ is constant on the set of all cones, and we will write $\Lambdaa {\text{cone}}$.

For \ssts, our main technical result is the following theorem.

\thm\label{thmMAINS} 
Let $\mcals$ be a set of \ssts\ that contains all pyramids and bipyramids.  Let $\Lambda$ be a \scvf\ on $\mcals$.  Suppose that $\Lambda$ is \gld\ by a \svf\ $\phi$ on $\mcals$, and that $\Lambda$ is constant on the set of all pyramids and bipyramids.  Then for each $K \in \mcals$ and each $v \in \verts K$, we have
\begin{equation}\label{eqAAS}
\phiw vK = \frac 12 \Lambdaa {\text{pyramid}}\left[1 - \sum_{\sigma^2 \ni v} \angg v{\sigma^2}\right],
\end{equation}
where the summation is over all \tsimps\ of $K$ containing $v$.
\ethm

The proof of Theorem~\ref{thmMAINS} will be given in Section~\ref{secPROO}.  It is straightforward to see that Theorem~\ref{thmSURF} follows immediately from Theorem~\ref{thmMAINS}.  Moreover, we have the following corollaries to Theorem~\ref{thmMAINS}, which characterize the Euler characteristic on the set of \ssts.  The first of these corollaries follows immediately from Theorem~\ref{thmMAINS} because of the \dgbt\ for the \cad, and the second corollary follows from the first.

\coro\label{coroMAINT} 
Let $\mcals$ be a set of \ssts\ that contains all pyramids and bipyramids.  Let $\Lambda$ be a \scvf\ on $\mcals$.  Suppose that $\Lambda$ is \gld\ and is constant on the set of all pyramids and bipyramids.  Then $\Lambda$ equals the Euler characteristic multiplied by $\frac 12 \Lambdaa {\text{pyramid}}$.
\ecoro

\coro\label{coroMAINV} 
Let $\mcals$ be a set of \ssts\ that contains all pyramids and bipyramids.  The Euler characteristic is the unique \scvf\ on $\mcals$ that is \gld\ and has value $2$ on all pyramids and bipyramids.
\ecoro


We now turn to the analogs of the above results for more general sets of \simpts, starting with the following definition, which is based on \cite{FORM5}.


\defn
Let $\mcalt$ be a set of \simpts.  The set $\mcalt$ is \textbf{\scl} if for every $K \in \mcalt$, and every vertex $v$ of $K$, we have $\stark vK \in \mcalt$.
\edefn

Clearly, the set of all finite \simpts\ is \scl.

We can now state the analog for arbitrary \simpts\ of Theorem~\ref{thmMAINS}.

\thm\label{thmMAINA} 
Let $\mcalt$ be a \scl\ set of \simpts\ that contains all planar fans.  Let $\Lambda$ be a \scvf\ on $\mcalt$.  Suppose that $\Lambda$ is \gld\ by a \svf\ $\phi$ on $\mcalt$, that $\Lambda$ is constant on the set of all planar fans, and for each $n \ge 0$ the function $\Lambda$ is constant on the set of all \nflaps n in $\mcalt$.  Then for each $K \in \mcalt$ and each $v \in \verts K$, we have
\begin{equation}\label{eqAAA}
\begin{aligned}
\phiw vK &= \Lambdaa {\stark vK} - \frac 12 \sum_{w \in \linkk vK} \Lambdaa {\text{\flaplinkord wv}}\\
 &\qquad\qquad + \frac 12 \Lambdaa {\text{fan}} f_1(\linkk vK) - \Lambdaa {\text{fan}} \sum_{\sigma^2 \ni v} \angg v{\sigma^2},
\end{aligned}
\end{equation}
where the first summation is over all vertices $w$ of $\linkk vK$, and the second summation is over all \tsimps\ of $K$ containing $v$.
\ethm

One could view Theorem~\ref{thmMAINA} as the geometric, \nd{2}dimensional analog of the uniqueness stated in \cite{LEVI}*{Theorem~B}.  The following result is an immediate consequence of Theorem~\ref{thmMAINA}

\coro\label{coroMAINE} 
Let $\mcalt$ be a \scl\ set of \simpts\ that contains all disks.  Let $\Lambda$ be a \scvf\ on $\mcalt$.  Suppose that $\Lambda$ is \gld\, that $\Lambda$ is constant on the set of all planar fans, and for each $n \ge 0$ the function $\Lambda$ is constant on the set of all \nflaps n in $\mcalt$.  Then there is a unique \svf\ $\phi$ on $\mcalt$ such that $\Lambda$ is \gld\ by $\phi$.
\ecoro

We note that Theorem~\ref{thmMAINA} implies not only that any \scvf\ $\Lambda$ that is \gld, and is constant on the set of all planar fans and is constant on the set of all \nflaps n\ for each $n \ge 0$, is \gld\ by a unique \svf\ $\phi$, but that such $\phi$ necessarily has the form of an angle defect, in that the first three terms in the right hand side of Equation~\ref{eqAAA} depends only upon $\linkk vK$ up to homeomorphism, and hence the right hand side of Equation~\ref{eqAAA} has the form of a measure of flatness (which is $1$ in the case of \ssts, and in general depends only upon the topology of a neighborhood of $v$) minus the sum of the angles at $v$ (once the term $\Lambdaa {\text{fan}}$ has been factored out).

It is straightforward to see that Theorem~\ref{thmSTCU} follows immediately from Corollary~\ref{coroMAINE}.  

The following two corollaries to Theorem~\ref{thmMAINA} characterize the Euler characteristic on the set of \ssts.  The first of these corollaries will be proved in  Section~\ref{secPROO}, and the second corollary follows from the first.

\coro\label{coroMAINC} 
Let $\mcalt$ be a \scl\ set of \simpts\ that contains all planar fans.  Let $\Lambda$ be a \scvf\ on $\mcalt$.  Suppose that $\Lambda$ is \gld, and that $\Lambda$ is constant on the set of all cones in $\mcalt$.  Then $\Lambda$ equals the Euler characteristic multiplied by $\Lambdaa {\text{cone}}$.
\ecoro

\coro\label{coroMAINCV} 
Let $\mcalt$ be a \scl\ set of \simpts\ that contains all planar fans.  The Euler characteristic is the unique \scvf\ on $\mcalt$ that is \gld, and has value $1$ on all cones in $\mcalt$.
\ecoro

\section{Proofs}
\label{secPROO}

We will prove Theorem~\ref{thmMAINA}, Corollary~\ref{coroMAINC} and Theorem~\ref{thmMAINS}.

\demoname{Proof of Theorem~\ref{thmMAINA}}
We start with the following preliminary observation.  Let $K$ and $L$ be \simpts, let $v \in \verts K$ and let $u \in \verts L$.  Suppose that $\linkk vK$ and $\linkk uL$ are both polygonal arcs, and that the sum of the angles at $v$ equals the sum of the angles at $u$.  Clearly there is a subdivision $K'$ of $K$ and a subdivision $L'$ of $L$ such that $\stark v{K'}$ and $\stark u{L'}$ are simplicially isometric.  It then follows from the invariance of $\phi$ under subdivision and under simplicial isometries of stars that $\phiw vK = \phiw v{K'} = \phiw u{L'} = \phiw uL$.

Our proof has a number of steps, first looking at some special cases, and then proving the result in general in the last step.  In each step, we will let $K$ be a \simpt, and we will let $v \in \verts K$ subject to certain stated conditions; we will then find a formula for $\phiw vK$ in the given case.  The arguments in many of the steps are similar to each other, and we will omit some of the details.

\stp\label{stp1c}
Suppose that $\linkk vK$ is a polygonal arc, and that the sum of the angles at $v$ is a number $\epsilon$ that has the form $\epsilon = \frac {n - 2}{2n}$ for some $n \in \nn$ such that $n \ge 3$.  We will show that $\phiw vK = \Lambdaa {\text{fan}}\left[\frac 12 - \epsilon\right]$.  (It can be verified that this last equation is a special case of Equation~\ref{eqAAA}, though we will not need this fact, and will not give the details.)

Let $S$ be a disk in $\rrr{2}$, the boundary of which is a regular polygon with $n$ vertices, say $a_1, \ldots, a_n$.  The angle at each $a_i$ equals $\frac {n - 2}{2n}$ (recall that we have normalized angles so that a complete circle has angle $1$).  The disk $S$ can be triangulated as a planar fan.  Then the link of each $a_i$ is a polygonal arc, and the sum of the angles at each $a_i$ is $\frac {n - 2}{2n}$.  By our preliminary observation, we know that all the $\phiw {a_i}S$ are equal to each other and to $\phiw vK$.  Applying the \dgbt\ to the disk $S$, we deduce that $\sum_{i=1}^n \phiw {a_i}S = \Lambdaa S$, and hence $n \phiw vK = \Lambdaa S$, which implies $\phiw vK = \frac 1n \Lambdaa {\text{fan}}$.  However, because $\epsilon = \frac {n - 2}{2n}$, we have $n = \frac 1{\frac 12 - \epsilon}$, and we deduce that $\phiw vK = \Lambdaa {\text{fan}}\left[\frac 12 - \epsilon\right]$.
\estp

\stp\label{stp1d}
Suppose that $\linkk vK$ is a polygonal arc, and that the sum of the angles at $v$ is a rational number $\delta$ such that $0 < \delta < \frac 12$.  We will show that $\phiw vK = \Lambdaa {\text{fan}}\left[\frac 12 - \delta\right]$.

Let $\delta = \frac pq$ for some $p, q \in \nn$.  Because $\delta < \frac 12$, then $q \ge 2$.  Let $D$ be a convex polygonal disk in $\rrr{2}$ with $p + 3$ vertices, labeled $s, x, t, a_1, \ldots, a_{p}$ in order around the boundary of $D$, such that the angles at $s$ and $t$ are $\frac 14$, the angle at $x$ is $\delta$, and the angle at each $a_i$ is $\frac {q - 2}{2q}$.  That such a convex polygon exists is due to the fact that all the angles are less than $\frac 12$, and the sum of the exterior angles is precisely $1$ (as can easily be verified).  The disk $D$ can be triangulated as a planar fan.  Observe that $\frac 14 = \frac {4 - 2}2 \cdot 4$, and hence all of the angles in $D$ other than the angle at $x$ satisfy the hypothesis of Step~\ref{stp1c}.  Applying the \dgbt\ to the disk $D$, and solving for $\phiw xD$, we deduce that 
$$
\phiw sD + \phiw xD + \phiw tD + \sum_{i=1}^{p} \phiw {a_i}D = \Lambdaa D.
$$
By Step~\ref{stp1c} we conclude that
\begin{align*}
&\phiw xD = \Lambdaa D - \phiw sD - \phiw tD - \sum_{i=1}^{p} \phiw {a_i}D \\
 &\ \ \ = \Lambdaa {\text{fan}} - \Lambdaa {\text{fan}}\left[\frac 12 - \frac 14\right]  - \Lambdaa {\text{fan}}\left[\frac 12 - \frac 14\right] - \sum_{i=1}^{p} \Lambdaa {\text{fan}}\left[\frac 12 - \frac {q - 2}{2q}\right]\\
 &\ \ \ = \Lambdaa {\text{fan}}\left[\frac 12 - \frac pq\right] = \Lambdaa {\text{fan}}\left[\frac 12 - \delta\right].
\end{align*}
Hence $\phiw vK = \Lambdaa {\text{fan}}\left[\frac 12 - \delta\right]$.
\estp

\stp\label{stp2}
Suppose that $\linkk vK$ is a polygonal arc, and that the sum of the angles at $v$ is a real number $\gamma$ such that $0 < \gamma < \frac 12$.  We will show that $\phiw vK = \Lambdaa {\text{fan}}\left[\frac 12 - \gamma\right]$.

There exists a sequence of positive rational numbers  $\{\delta_n\}_{n=1}^\infty$ such that $0 < \delta_n < \frac 12$ for all $n$, and $\lim_{n \to \infty} \delta_n = \gamma$.  Let $T = \langle x, y, z \rangle$ be a triangle in $\rrr{2}$ such that the angle at $x$ is $\gamma$.  Clearly there is a sequence of triangles $\{T_n\}_{n=1}^\infty$ in $\rrr{2}$, where for each $n$ we have $T_n = \langle x_n, y_n, z_n \rangle$ with the angle at $x_n$ equal to $\delta_n$, and such that $\lim_{n\to\infty} T_n = T$, with $\lim_{n\to\infty} x_n = x$.  By Step~\ref{stp1d} we know that $\phiw {x_n}{T_n} = \Lambdaa {\text{fan}}\left[\frac 12 - \delta_n\right]$ for all $n$.    By the continuity of $\phi$, we know that $\phiw xT = \lim_{n\to\infty} \phiw {x_n}{T_n} = \lim_{n\to\infty} \Lambdaa {\text{fan}}\left[\frac 12 - \delta_n\right] = \Lambdaa {\text{fan}}\left[\frac 12 - \gamma\right]$.  Hence $\phiw vK = \Lambdaa {\text{fan}}\left[\frac 12 - \gamma\right]$.
\estp

\stp\label{stp3}
Suppose that $\linkk vK$ is a polygonal arc, and that the sum of the angles at $v$ is a real number $\beta$ such that $0 < \beta < 1$.  We will show that $\phiw vK = \Lambdaa {\text{fan}}\left[\frac 12 - \beta\right]$.

Let $Q$ be a quadrilateral in $\rrr{2}$ with vertices $e, b_1, b_2, b_3$, such that the angle at $e$ is $\beta$, and the other three angles are less than $\frac 12$.  The quadrilateral $Q$ can be triangulated as a planar fan.  The desired result can be obtained by applying the \dgbt\ to the disk $Q$, solving for $\phiw eQ$, and then using Step~\ref{stp2} and the fact that the sum of the angles in a quadrilateral is $1$.
\estp

\stp\label{stp4}
Suppose that $K$ is an \nflap n\ for some $n \ge 0$, and that $v$ is one of the \everts\ of $K$.  We will show that
\begin{equation}\label{eqAAC}
\phiw vK = \frac 12 \Lambdaa {\text{\nflap n}} - \Lambdaa {\text{fan}} \sum_{\sigma^2 \ni v} \angg v{\sigma^2},
\end{equation}
where the  summation is over all \tsimps\ of $K$ containing $v$.  (Again, it can be verified that Equation~\ref{eqAAC} is a special case of Equation~\ref{eqAAA}, though we will not need this fact.)

Let $w$ be the other \evert\ of $K$.  We have four subcases.

\setstepn{0}

\stpn{Subcase}\label{case11}
Suppose that $n = 0$.  Then $K$ consists of a single edge $\ip vw$ together with its vertices.  Clearly $v$ and $w$ have simplicially isometric stars, and hence $\phiw vK = \phiw wK$.  Applying the \dgbt\ to $\phi$, we deduce that $\phiw vK = \frac 12 \Lambdaa {\ip vw}$.  Because $\ip vw$ is a \nflap 0, and $v$ is contained in no \tsimps, clearly Equation~\ref{eqAAC} holds in this case.
\estpn

\stpn{Subcase}\label{case12}
Suppose that $n = 1$.  Then $K$ is a single \tsimp.  Let $\omega$ be the angle at $v$ in $K$.  We deduce from Step~\ref{stp2} that $\phiw vK = \Lambdaa {\text{fan}}\left[\frac 12 - \omega\right]$.  This last equation is a special case of Equation~\ref{eqAAC}, using the fact that a \nflap 1\ is a planar fan.
\estpn

\stpn{Subcase}\label{case13}
Suppose that $n \ge 2$.  Assume that $\angg v{\sigma^2} < \frac 14$ for all \tsimps\ $\sigma^2$ of $K$ containing $v$.

First, observe that because $K$ is an \nflap n, we can find an embedding of $K$ in $\rrr{3}$ that is simplicially isometric with $K$.  Hence, by the invariance of $\phi$ under simplicial isometries of stars, we may assume without loss of generality that $K$ is in $\rrr{3}$.

Choose a plane $\Pi$ in $\rrr{3}$ that intersects the relative interior of the edge $\ip vw$ and is perpendicular to it, and such $v$ is on one side of $\Pi$ and all the other vertices of $K$ are on the other side.  Let $V$ be the intersection of $|K|$ with the closed half-space in $\rrr{3}$ that has $\Pi$ as its boundary and contains $v$.  Let $V'$ be the the result of reflecting $V$ in $\Pi$, and let $Y = V \cup V'$.  The fact that $Y$ is an \nflap n\ follows from the hypothesis concerning the angles in the \tsimps\ of $K$ that contain $v$, and the choice of $\Pi$.  Let $x$ denote the vertex of $Y$ that is the mirror image of $v$, let $\tau_1, \ldots, \tau_n$ denote the $n$ \tsimps\ in $Y$, and for each $i \in \{1, \ldots, n\}$ let $d_i$ denote the vertex in $\tau_i$ that is not $v$ or $x$.  

By Step~\ref{stp3}, and making use of the symmetry of $Y$, we see that $\phiw {d_i}Y = \Lambdaa {\text{fan}}\left[\frac 12 - \angg {d_i}Y\right] = 2\Lambdaa {\text{fan}}\angg v{\tau_i}$ for each $i \in \{1, \ldots, n\}$.  By the invariance of $\phi$ under subdivision and under simplicial isometries of stars we see that $\phiw vK = \phiw vY = \phiw xY$.  The \dgbt\ applied to $Y$ then implies that 
$$
2\phiw vY + \sum_{i=1}^n \phiw {d_i}Y = \Lambdaa {Y},
$$
and hence that
\begin{align*}
\phiw vK &= \phiw vY = \frac 12\Bigl\{\Lambdaa {Y} - \sum_{i=1}^n \phiw {d_i}Y \Bigr\}\\
 &= \frac 12 \Lambdaa {\text{\nflap n}} - \Lambdaa {\text{fan}} \sum_{i=1}^n \angg v{\tau_i},
\end{align*}
which is equivalent to Equation~\ref{eqAAC}.
\estpn

\stpn{Subcase}\label{case14}
Suppose that $n \ge 2$, but we make no assumptions regarding the angles in $K$.  (Note that the argument used in Subcase~\ref{case13} will not work in the general case, because $Y$ as constructed would not always be an \nflap n, though it would have underlying space homeomorphic to one.)  As in Subcase~\ref{case13}, we may assume without loss of generality that $K$ is in $\rrr{3}$.

We now modify $K$ as follows.  Let $\zeta_1, \ldots, \zeta_n$ denote the $n$ \tsimps\ in $K$, and for each $i \in \{1, \ldots, n\}$ let $e_i$ denote the vertex in $\zeta_i$ that is not $v$ or $w$.  Let $i \in \{1, \ldots, n\}$.  If $\angg w{\zeta_i} < \frac 14$, then we leave $\zeta_i$ unchanged.  If $\angg w{\zeta_i} \ge \frac 14$, then we modify $\zeta_i$ by moving $e_i$ along the line segment $\ip v{e_i}$ toward $v$ until $\angg w{\zeta_i}$ becomes less than $\frac 14$.  This modification will not change $\angg v{\zeta_i}$.  After we modify all triangles as needed, we call the new \nflap n\ $Z$.  Similarly to previous arguments, we note that $\phiw vK = \phiw vZ$.  Observe also that the vertex $w$ in $Z$ satisfies the hypotheses of Subcase~\ref{case13}, and so $\phiw wZ$ satisfies Equation~\ref{eqAAC}.  

By Step~\ref{stp3}, for each $i \in \{1, \ldots, n\}$ we have 
$$
\phiw {e_i}Z = \Lambdaa {\text{fan}}\left[\frac 12 - \angg {e_i}Z\right] = \Lambdaa {\text{fan}}\left[\angg v{\zeta_i} + \angg w{\zeta_i}\right].
$$
The \dgbt\ applied to $Z$, and then solving for $\phiw vZ$, we see that
\begin{align*}
\phiw vK &= \phiw vZ = \Lambdaa {Z} - \phiw wZ - \sum_{i=1}^n \phiw {e_i}Z\\
 &= \Lambdaa {\text{\nflap n}} - \Bigl\{\frac 12 \Lambdaa {\text{\nflap n}} - \Lambdaa {\text{fan}} \sum_{i=1}^n \angg w{\zeta_i}\Bigr\}\\
 &\qquad\qquad - \Lambdaa {\text{fan}} \sum_{i=1}^n \left[\angg v{\zeta_i} + \angg w{\zeta_i}\right],
\end{align*}
which implies Equation~\ref{eqAAC}.
\estpn
\estp

\stp\label{stp5}
We now prove our theorem.  Suppose that $K$ is a \simpt, and that $v \in \verts K$.

Let $N = \stark vK$.  Clearly $\linkk vK = \linkk vN$, and by the invariance of $\phi$ under simplicial isometries of stars we know that $\phiw vK = \phiw vN$.  Hence we can rewrite Equation~\ref{eqAAA} as
\begin{equation}\label{eqAAAX}
\begin{aligned}
\phiw vN &= \Lambdaa {N} - \frac 12 \sum_{w \in \linkk vN} \Lambdaa {\text{\flaplinkord wv}} + \frac 12 \Lambdaa {\text{fan}} f_1(\linkk vN)\\
 &\qquad\qquad - \Lambdaa {\text{fan}} \sum_{\sigma^2 \ni v} \angg v{\sigma^2},
\end{aligned}
\end{equation}
where the first summation is over all vertices $w$ of $\linkk vN$, and the second summation is over all \tsimps\ of $\stark vN$ containing $v$.  We will prove Equation~\ref{eqAAAX}. 

If $\linkk vN = \emptyset$, then $N = \{v\}$, and Equation~\ref{eqAAAX} is trivially true, because the \dgbt\ applied to $N$ yields $\phiw vN = \Lambdaa {N}$, and $\linkk vN = \emptyset$ implies that all the terms in the right hand side of Equation~\ref{eqAAAX} other than $\Lambdaa {N}$ are zero.  Hence we suppose that $\linkk vN \ne \emptyset$.

Let $r$ be a vertex in $\linkk vN$.  Let $M = \stark rN$.  As remarked after the definition of \nflaps n, we know that $M$ is a \flaplinkord rv\ with \everts\ $v$ and $r$.  By the invariance of $\phi$ under simplicial isometries of stars, we know that $\phiw rN = \phiw rM$.  It then follows from Step~\ref{stp4} that  
\begin{equation}\label{eqAACX}
\phiw rN = \phiw rM = \frac 12 \Lambdaa {\text{\flaplinkord rv}} - \Lambdaa {\text{fan}} \sum_{\sigma^2 \ni r} \angg r{\sigma^2},
\end{equation}
where the summation is over all \tsimps\ of $M$ containing $r$.

Recall that all the vertices of $N$ other than $v$ are in $\linkk vN$.  Applying the \dgbt\ to $N$, and then solving for $\phiw vN$, we see that
\begin{align*}
&\phiw vN = \Lambdaa {N} - \sum_{w \in \linkk vN} \phiw wN\\
 &\qquad= \Lambdaa {N} - \sum_{w \in \linkk vN} \Bigl\{ \frac 12 \Lambdaa {\text{\flaplinkord wv}} - \Lambdaa {\text{fan}} \sum_{\sigma^2 \ni w} \angg w{\sigma^2}\Bigr\}\\
\intertext{\hfill by Equation~\ref{eqAACX}}
 &\qquad= \Lambdaa {N} - \frac 12 \sum_{w \in \linkk vN} \Lambdaa {\text{\flaplinkord wv}} + \!\!\!\!\sum_{w \in \linkk vN} \Lambdaa {\text{fan}} \sum_{\sigma^2 \ni w} \angg w{\sigma^2}\\
 &\qquad= \Lambdaa {N} - \frac 12 \sum_{w \in \linkk vN} \Lambdaa {\text{\flaplinkord wv}} + \!\!\!\!\sum_{\sigma^2 \in \stark vN} \Lambdaa {\text{fan}} \sum_{\substack {y \in \sigma^2 \\ y \ne v}} \angg y{\sigma^2}\\
\intertext{\hfill because every \nd{2}simplex in $N$ contains $v$}
 &\qquad= \Lambdaa {N} - \frac 12 \sum_{w \in \linkk vN} \Lambdaa {\text{\flaplinkord wv}} + \!\!\!\!\sum_{\sigma^2 \in \stark vN} \Lambdaa {\text{fan}} \left[\frac 12 - \angg v{\sigma^2}\right]\\
\intertext{\hfill because the sum of the angles in a triangle is $\frac 12$}
 &\qquad= \Lambdaa {N} - \frac 12 \sum_{w \in \linkk vN} \Lambdaa {\text{\flaplinkord wv}}
+ \frac 12 \Lambdaa {\text{fan}} f_1(\linkk vN)\\
 &\qquad\qquad\qquad - \Lambdaa {\text{fan}} \sum_{\sigma^2 \ni v} \angg v{\sigma^2},
\end{align*}
where the last equation holds because $f_2(\stark vN) = f_1(\linkk vN)$.  Hence Equation~\ref{eqAAAX} holds.
\estp
\edemoname 

For our next proof, we will need the following notation.  For any simplex $\eta$ of dimension $0$, $1$ or $2$ in Euclidean space, and any vertex $v$ of $\eta$, we let $\angge v\eta$ denote the exterior angle of $\eta$ at $v$.  If 
$\eta$ is a \tsimp, then $\angge v\eta = \frac 12 - \angg v\eta$; if $\eta$ is a \nd{1}simplex, then $\angge v\eta = \frac 12$; if $\eta$ is a \nd{0}simplex (so that $\sigma = v$), then $\angge v\eta = 1$.

\demoname{Proof of Corollary~\ref{coroMAINC}}
Note that every planar fan, and every star of every vertex of a simplicial complex in $\mcalt$, and in particular every \flaplinkord wv\ for appropriate vertices $v$ and $w$, are cones.

Let $K \in \mcalt$, and let $v \in \verts K$.  It now follows from Equation~\ref{eqAAA} that
$$
\phiw vK = \Lambdaa {\text{cone}}\Bigl\{1 - \frac 12 f_0(\linkk vK) + \frac 12 f_1(\linkk vK) - \sum_{\sigma^2 \ni v} \angg v{\sigma^2}\Bigr\}.
$$

As given in a number of sources, for example \cite{BA1}, we know that the \stcurv\ of $K$ at $v$ is given by
$$
\kbk vK = \sum_{i=0}^2 (-1)^i \sum_{\eta^i \ni v} \angge v{\eta^i},
$$
where the inner summation is over all \nd{i}simplices of $K$ containing $v$.  We then compute
\begin{align*}
\allowdisplaybreaks
\kbk vK &= \sum_{i=0}^2 (-1)^i \sum_{\eta^i \ni v} \angge v{\eta^i}\\
 &= \angge vv - \sum_{\eta^1 \ni v} \angge v{\eta^1} + \sum_{\eta^2 \ni v} \angge v{\eta^2}\\
 &= 1 - \sum_{\eta^1 \ni v} \frac 12 + \sum_{\eta^2 \ni v} \left[\frac 12 - \angg v{\eta^2}\right]\\
 &=1 - \frac 12 f_0(\linkk vK) + \sum_{\eta^2 \ni v} \frac 12 - \sum_{\eta^2 \ni v} \angg v{\eta^2}\\
 &=1 - \frac 12 f_0(\linkk vK) + \frac 12 f_1(\linkk vK) - \sum_{\eta^2 \ni v} \angg v{\eta^2}.
\end{align*}
It follows that $\phiw vK = \Lambdaa {\text{cone}}\kbk vK$.  

Because \stcurv\ satisfies the \dgbt\ with respect to the Euler characteristic, and because $\phi$ satisifies the \dgbt\ with respect to $\Lambda$, it follows that $\Lambda$ equals the Euler characteristic multiplied by the constant $\Lambdaa {\text{cone}}$.
\edemoname

We now turn to the proof of Theorem~\ref{thmMAINS}.  The proof of this theorem is not identical to the proof of Theorem~\ref{thmMAINA}, because in the former theorem we assume that $\phi$ is defined only on \ssts, whereas in the proof of the latter theorem we make use of various simplicial complexes that are not \ssts.  We will show, however, how the proof of Theorem~\ref{thmMAINA} can be modified to work in the case of \ssts.  (It would be easier to prove Theorem~\ref{thmMAINS} if we allowed non-embedded \ssts, but that would add unnecessarily to the hypotheses of the theorem, so we will not do so.)

\demoname{Proof of Theorem~\ref{thmMAINS}}
This proof has a number of steps, most of which are similar to some of the steps of the proof of Theorem~\ref{thmMAINA}.  We start with some observations.

\stpnn{(a)}
Because all simplicial complexes under consideration are \ssts, we know that the link of every vertex is a polygonal circle.
\estpnn

\stpnn{(b)}
Let $K, L \in \mcals$, let $v \in \verts K$ and let $w \in \verts L$.  Suppose that the sum of the angles at $v$ equals the sum of the angles at $w$.  Clearly there is a subdivision $K'$ of $K$ and a subdivision $L'$ of $L$ such that $\stark v{K'}$ and $\stark w{L'}$ are simplicially isometric.  It then follows from the invariance of $\phi$ under subdivision and under simplicial isometries of stars that $\phiw vK = \phiw wL$.  

\begin{figure}[ht]
\hfil\includegraphics{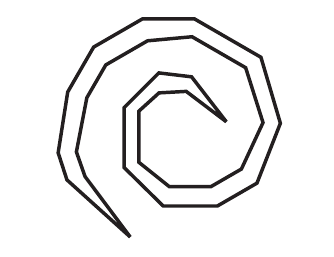}\hfil\hfil
\caption{}\label{figAXIa}
\end{figure}

Let $\omega \in (0, \infty)$.  Then it is possible to draw a polygonal spiral ribbon $R$ in $\rrr{2}$, as in Figure~\ref{figAXIa}, so that an appropriately chosen bipyramid $B$ on $R$ has angle sum $\omega$ at each of the apices, and has angle sum less than $1$ at all other vertices (this latter condition will be used later in the proof).  Observe that $B \in \mcals$, because $\mcals$ contains all bipyramids.

We can therefore think of $\phi$ as a function $(0, \infty) \to \rr$, where for each $\alpha \in (0, \infty)$  we define $\phia \alpha$ by $\phia \alpha = \phiw vK$ for any $K \in \mcals$ that has a vertex $v$ for which the sum of the angles at $v$ is $\alpha$.
\estpnn

\stpnn{(c)}
The continuity of $\phi$ as originally assumed implies that $\phi$ is continuous when thought of as a function $(0, \infty) \to \rr$, as described in Observation~(b).  
\estpnn

Given the above observations, in order to prove the theorem as originally stated it suffices to show that
\begin{equation}\label{eqAASW}
\phia \omega = \frac 12 \Lambdaa {\text{pyramid}}\left[1 - \omega\right],
\end{equation}
for all $\omega \in (0, \infty)$.

For the rest of this proof, let $\omega \in (0, \infty)$.  We have a number of cases.

\setstepn{0}

\stpn{Case}\label{stp1ax}
Suppose that $\omega = 1$.    

Let $T$ be a triangle in $\rrr{3}$, with vertices $d_1, d_2, d_3$.  Let $\Delta$ be a pyramid on $T$ with apex $b$.  Let $\delta_i$ denote the sum of the angles in $\Delta$ at $d_i$, and let $\beta$ denote the sum of the angles in $\Delta$ at $b$.  Let $\{L_k\}_{k=1}^\infty$ in $\rrr{3}$ be a sequence of bipyramids on $T$, where all the $L_k$ have $b$ as one of their cone apices, where the other apex in $L_k$ is denoted $e_k$ for each $k$, and where the sequence $\{e_k\}_{k=1}^\infty$ converges to the centroid of $T$.  For each $k$, let $\delta_i^k$ denote the sum of the angles in $L_k$ at $d_i$, and let $\gamma^k$ denote the sum of the angles in $L_k$ at $e_k$; observe that the sum of the angles in each $L_k$ at $b$ is $\delta$.  Clearly, we see that $\lim_{k\to\infty} \delta_i^k = \delta_i$ for all $i$, and $\lim_{k\to\infty} \gamma^k = 1 = \omega$.

Let $k \in \nn$.  By appyling the \dgbt\ to $L_k$ we obtain
$$
\phiw b{L_k} + \sum_{i=1}^3 \phiw {d_i}{L_k} + \phiw {e_k}{L_k} = \Lambdaa {L_k}.
$$
Because $\Lambda$ is constant on the set of all pyramids and bipyramids, we know that $\Lambdaa {L_k} = \Lambdaa \Delta$ for all $k$.  It then follows that
$$
\phia \beta + \sum_{i=1}^3 \phia {\delta_i^k} + \phia {\gamma^k} = \Lambdaa \Delta.
$$
Taking the limit as $k \to \infty$, and using the continuity of $\phi$ as stated in Observation~(c), we see that 
\begin{equation}\label{eqABA}
\phia \beta + \sum_{i=1}^3 \phia {\delta_i} + \phia \omega = \Lambdaa \Delta.
\end{equation}
On the other hand, appyling the \dgbt\ to $\Delta$ yields
$$
\phiw b{\Delta} + \sum_{i=1}^3 \phiw {d_i}{\Delta} = \Lambdaa {\Delta},
$$
which implies
\begin{equation}\label{eqABB}
\phia \beta + \sum_{i=1}^3 \phia {\delta_i} = \Lambdaa \Delta.
\end{equation}
Comparing Equation~\ref{eqABA} with Equation~\ref{eqABB} shows that $\phia \omega = 0$, which is equivalent to Equation~\ref{eqAASW} in the present case.
\estpn

\stpn{Case}\label{stp1cx}
Suppose that $\omega = \frac {n - 2}n$ for some $n \in \nn$ such that $n \ge 3$.  

Let $S$ and $a_1, \ldots, a_n$ be as in Step~\ref{stp1c} of the proof of Theorem~\ref{thmMAINA}.  Let $\{C_k\}_{k=1}^\infty$ in $\rrr{3}$ be a sequence of pyramids on $S$, where the apex in $C_k$ is denoted $x_k$ for each $k$, where the \tsimps\ of $C_k$ containing $x_k$ are all congruent isosceles triangles, and where the sequence $\{x_k\}_{k=1}^\infty$ converges to the centroid of $S$.  For each $k$, let $\omega_i^k$ denote the sum of the angles in $C_k$ at $a_i$, and let $\beta^k$ denote the sum of the angles in $C_k$ at $x_k$.  Clearly $\lim_{k\to\infty} \omega_i^k = \omega$ for all $i$, and $\lim_{k\to\infty} \beta^k = 1$.

We now proceed similarly to Case~\ref{stp1ax}.  Let $k \in \nn$.  By appyling the \dgbt\ to $C_k$ we obtain
$$
\sum_{i=1}^n \phia {\omega_i^k} + \phia {\beta^k} = \Lambdaa {C_k}.
$$
By construction we know that the $\omega_i^k$ are all equal to each other for all $i$, and hence 
$$
\phia {\omega_1^k}= \frac 1n \left[\Lambdaa {\text{pyramid}} - \phia {\beta^k}\right].
$$
Taking the limit as $k \to \infty$, and using the continuity of $\phi$, as well as Case~\ref{stp1ax}, we see that 
\begin{equation}\label{eqE2W}
\phia {\omega}= \frac 1n \left[\Lambdaa {\text{pyramid}} - \phia {1}\right] = \frac 1n \Lambdaa {\text{pyramid}}.
\end{equation}
Because $\omega = \frac {n - 2}n$, we have $n = \frac 2{1 - \omega}$.  Substituting this formula for $n$ into Equation~\ref{eqE2W}, we see that Equation~\ref{eqAASW} holds in this case.
\estpn

\stpn{Case}\label{stp1dx}
Suppose that $\omega$ is a rational number such that $0 < \omega < 1$.  

There are $p, q \in \nn$ such that $\omega = \frac {2}q$.  Because $\omega < 1$, then $q \ge 2$.  The argument used to show that $\phia \omega$ is given by Equation~\ref{eqAASW} is similar to the argument in Case~\ref{stp1cx}, except that we take pyramids on the polygon $D$ in Step~\ref{stp1d} of the proof of Theorem~\ref{thmMAINA}; we omit the details.
\estpn

\stpn{Case}\label{stp2x}
Suppose that $\omega \in (0, 1)$.  

Equation~\ref{eqAASW} follows immediately from Case~\ref{stp1dx} and the continuity of $\phi$.
\estpn

\stpn{Case}\label{stp3x}
Suppose that $\omega \in [1, \infty)$.  

As remarked in Observation~(b), it is possible to draw a polygonal spiral ribbon $R$ in $\rrr{2}$ so that an appropriately chosen bipyramid $B$ on $R$ has angle sum $\omega$ at each of the apices, and has angle sum less than $1$ at all other vertices.  Suppose that the vertices of $R$ are denoted $b_1, \ldots, b_m$, and the apices of $B$ are denoted $x$ and $y$.  Suppose further that for each $i$, the sum of the angles in $B$ at $b_i$ is $\beta_i$.  By Case~\ref{stp2x} we know that $\phia {\beta_i} = \frac 12 \Lambdaa {\text{pyramid}}\left[1 - \beta_i\right]$ for each $i$.

Clearly the \cad\ at $b_i$ is $1 - \beta_i$, and the \cad\ at each of $x$ and $y$ is $1 - \omega$.  Using the \dgbt\ for the \cad\ applied to $B$, we have
$$
2 = \chi (B) = \sum_{i=1}^m \left[1 - \beta_i\right] + 2\left[1 - \omega\right],
$$
which implies that 
$$
\omega = \frac 12\sum_{i=1}^m \left[1 - \beta_i\right].
$$

Applying the \dgbt\ to $B$ yields
$$
\Lambdaa B = \sum_{i=1}^m \phiw {b_i}B + \phiw xB + \phiw yB,
$$
and using arguments similar to those used previously in this proof, we deduce that
\begin{align*}
\phia \omega &= \frac 12 \Bigl\{\Lambdaa {\text{pyramid}} - \sum_{i=1}^m \phia {\beta_i}\Bigr\}\\
 &= \frac 12 \Lambdaa {\text{pyramid}} \Bigl\{1 - \frac 12\sum_{i=1}^m \left[1 - \beta_i\right]\Bigr\} = \frac 12 \Lambdaa {\text{pyramid}}\left[1 - \omega\right].
\end{align*}
\estpn
\edemoname

\end{document}